\documentclass[12pt]{amsart}
\usepackage{amscd,amssymb}
\usepackage[arrow,matrix]{xy}

\theoremstyle{plain}
\newtheorem{thm}[subsection]{Theorem}

\newtheorem{prop}[subsection]{Proposition}
\newtheorem{cor}[subsection]{Corollary}

\theoremstyle{definition}
\newtheorem{rk}[subsection]{Remark}

\newtheorem{ex}[subsection]{Example}

\numberwithin{equation}{section}
\newcommand{\OO}{{\mathcal O}}

\renewcommand{\SS}{{\mathcal S}}

\newcommand{\Z}{\mathbb{Z}}

\newcommand{\C}{\mathbb{C}}

\newcommand{\PP}{\mathbb{P}}

\begin{document}

\title [On the connectivity of some complete intersections ]
{On the connectivity of some complete intersections     }

\author[Alexandru Dimca]{Alexandru Dimca }
\address{  Laboratoire J.A. Dieudonn\'e, UMR du CNRS 6621,
                 Universit\'e de Nice-Sophia-Antipolis,
                 Parc Valrose,
                 06108 Nice Cedex 02,
                 FRANCE.}
\email
{dimca@math.unice.fr}

\subjclass[2000]{Primary 14F25, 14F35; Secondary 32S20.}

\keywords{rational connectivity, complete intersection}

\begin{abstract}
We show that the complement of a degree $d$ hypersurface in a projective complete intersection,
whose defining equations have degrees strictly larger than $d$, has a rational connectivity higher than expected. The key new feature is that a positivity result replaces the usual transversality conditions needed
to get such connectivity results.

\end{abstract}

\maketitle

\section{Introduction} \label{sec:intro}

Let $V$ be an $n$-dimensional complex complete intersection in $\PP^N$, defined by the equations $f_1=f_2= \cdots = f_k=0$. Here $N=n+k$ and each polynomial $f_j \in \C[X_0,  \cdots , X_N]$ is homogeneous of some degree $d_j$. We assume that the singular set $V_{sing}$ is finite.
Let $H \subset \PP^N$ be a reduced hypersurface given by $h=0$, for a 
homogeneous polynomial $h$ of degree $d$.
We define the tangency set of $V$ and $H$ to be the set
$$T(V,H)=\{x \in V_{reg} \cap  H_{reg} ~ | ~ T_xV \subset T_xH ~ \}$$
where $V_{reg} =V \setminus V_{sing}$ is the smooth part of $V$ and $ H_{reg} $ has a similar meaning. 
The main result of this note is the following.

\begin{thm} \label{thm1}
Assume that $V_{sing} \cap H =\emptyset$ and   $d < \min _j d_j$. Then the following hold.

\medskip

\noindent (i) The tangency set  $T(V,H)$ is finite. In particular, if in addition $ H_{sing} \cap V$ is a finite set, then the complete intersection $W=V \cap H$ has at most isolated singularities.

\medskip

\noindent (ii) Assume that the Milnor fiber $F_h=\{x \in \C^{N+1}~ | ~ h(x)=1~\}$ is $s$-connected, for some $s \leq n-1$. Then
the Milnor fiber of the function germ $g:(CV,0) \to (\C,0)$, namely the  $n$-dimensional affine complete intersection
$$F=\{x \in \C^{N+1}~ | ~ f_1(x)=f_2(x)= \cdots = f_k(x)=0, ~ h(x)=1 ~\}$$
is also $s$-connected.  In particular, if $s=n-1$ e.g. if $\dim H_{sing}<k$, then $F$
is a bouquet of $n$-dimensional spheres.

\end{thm}

\begin{rk} \label{rk1} \rm

Without any assumption on the degrees, if $W$ is smooth, then the second claim above is classical, see Hamm \cite{H1} when $V$ is smooth or refer to Looijenga \cite{L}, p. 75 for the general case. Similarly, when  $W$ has at most isolated singularities, by taking a generic hyperplane section through the origin and using \cite{H2}, we get that $F$ is $(n-2)$-connected. The whole point of the paper is that the assumption $d < \min _j d_j$ replaces in a misterious way the transversality of $H$ with respect to $V$ and leads to stronger connectivity results than the general ones just described in this Remark.

\end{rk}

\begin{ex} \label{ex1} \rm
Assume that $k=1$ and $d=1$, i.e. $V$ is a  hypersurface having at most isolated singularities and $H$ is a hyperplane avoiding these singularities. Both claims are then already known,
see for instance \cite{D0}, p. 205   or \cite{CD} for the first claim and  \cite{Dtop} for the second (this case is in fact just a reformulation of A. N\'emethi results on quasi-tame polynomials in 
\cite{N1}, \cite{N2}). More generally, for a smooth complete intersection which is non-degenerate in the sense that $d_j>1$ for all $j=1,  \cdots, k$ and a hyperplane, the first claim is just Remark
7.5 in  \cite{FL}, see also  \cite{CD} where the claim (i) is obtained for $k=1$.

\end{ex}

\begin{ex} \label{ex2} \rm
Assume that $V$ and $H$ are smooth hypersurfaces of the same degree.
Then it is easy to see that the singular locus of $W$ can be arbirarily large. Indeed, for
$$f_1=X_0^2+X_1^2+\cdots +X_N^2, ~ ~ ~ h=2X_0^2+X_1^2+\cdots +X_N^2$$
$W$ is not even reduced, while for 
$$f_1=X_0^2+X_1^2+\cdots +X_N^2, ~ ~ ~ h=2X_0^2+2X_1^2+X_2^2+ \cdots +X_N^2$$
$W$ is reduced, with a codimension one singular set.

Moreover, even when $W$ has at most isolated singularities, this does not imply that $F$ 
is a bouquet of $n$-dimensional spheres. The following example was kindly provided by
 A. N\'emethi. When $N=5$, consider the equations
$$f_1=X_0^2+X_1^2+\cdots +X_5^2, ~ ~ ~ h=X_2^2+X_3^2+2X_4^2+2X_5^2.$$
Then $H_3(F)=\Z^3$ even though $n=4$ in this case. Of course, $F$ is 2-connected by our
Remark \ref{rk1}.
To obtain the same drop in connectivity in the case $d_1<d$, it is enough to take $d_1=1$
and note that $F$ is in this case just the complement $M=V \setminus W=\PP^n  \setminus W$ of a projective hypersurface having at most isolated singularities. Examples of such complements with $H_{n-1}(M) \ne 0$ can be easily obtained,
see for instance \cite{D1}, Chapter 6.

\end{ex}

Exactly as in the case $k=1$ treated in  \cite{Dtop} we get from Theorem \ref{thm1} the following. Recall the convention
$\dim \emptyset =-1$.

\begin{cor} \label{cor1}

Let $X$ be an $n$-dimensional complete intersection in $\C^N$ and let $V$ be the projective closure of $X$ in  $\PP^N$. If $H_{\infty}$ denotes the hyperplane at infinity $\PP^N \setminus \C^N$
and $\delta = \dim (V_{sing} \cap H_{\infty})$, then $X$ is $(n-\delta -2)$-connected.

\end{cor}

The following consequence is also obvious, compare to \cite{D1}, p.146.

\begin{cor} \label{cor2}
With the notation and assumptions from Theorem  \ref{thm1}, there is an $s$-equivalence
$\phi: M:=V \setminus W \to K(\mu_d,1)$, where $\mu_d$ is the cyclic group of $d$-roots of unity.
In particular, the reduced rational cohomology of $M$ vanishes in degrees $\leq s$.

\end{cor}

\section{Tangency sets} \label{s1}

In this section we prove the first claim in Theorem  \ref{thm1}. The method is already present in 
  \cite{CD}, but for the reader's convenience we recall it here. This idea plays also in key role in proving the
topological claim in Theorem  \ref{thm1} in the next section.

Assume that $\dim T(V,H)>0$. Then we can find, by taking repeated hyperplane sections if necessary, an irreducible algebraic set $C$ of dimension 1 such that $C \subset  T(V,H)$. For $x \in \C^{N+1}$, $x\ne 0$, we denote by $[x]$ the corresponding point in $\PP^N$. For $[x]\in C$ we have
\begin{equation} \label{eq1}
d_xh\in Span(d_xf_1, \cdots, d_xf_k).
\end{equation}
Moreover, if $[x_n]\in C$ is a sequence of points in $C$ converging to $[y]\in {\overline C}$, then
by continuity and using the assumption  $V_{sing} \cap H =\emptyset$, we get
\begin{equation} \label{eq2}
d_yh\in Span(d_yf_1, \cdots, d_yf_k).
\end{equation}
For $i=0, \cdots, N$, let $U_i=\{ [x]\in \PP^N ~ |~ x_i \ne 0 ~\}$. Let $D=\max _j d_j$ and define the following
homogeneous polynomials (all of degree $D$):
$h^i=X_i^{D-d}h(X)$ and $f_j^i=X_i^{D-d_j}f_j(X)$ for $j=1, \cdots, k.$ For $ [x] \in {\overline C}\cap U_i$, it follows that
$$d_xh^i=\sum_j\lambda _j^id_xf_j^i$$
where the complex numbers $\lambda _j^i$ depend only on the point $[x]$, and not on the chosen representative $x$. For
$ [x] \in {\overline C}\cap U_i \cap U_j$, we get
\begin{equation} \label{eq3}
\frac{\lambda _a^i}{x_i^{d_a-d}}=\frac{\lambda _a^j}{x_j^{d_a-d}}
\end{equation}
for all $a=1, \cdots, k.$ Therefore, the collection of regular functions $(\lambda _a^i)_i$ defines a regular section 
$\lambda_a$ of the line bundle $\OO (d-d_a)$ over the irreducible projective curve $ {\overline C}$.
Since $d-d_a<0$, it follows that this section is identically zero on  $ {\overline C}$. Applying this to all sections
$\lambda_a$, we get that $C \subset H_{sing} $, a contradiction. Indeed, by definition $C \subset T(V,H)$ and
$ T(V,H) \cap H_{sing} =\emptyset$.

\section{Homotopy types of Milnor fibers} \label{s2}

In this section we prove the second claim in Theorem  \ref{thm1}. The proof is divided in several steps,
so that the reader may easily follow the argument.

\subsection{Step 1: finding good equations for $V$} \label{s2a}

Recall the Bertini Theorem, saying that the general member of a linear system is smooth outside the base locus,
see for instance  \cite{GH}, p. 137. Assume that the degrees $d_j$ of the equations defining $V$ are ordered such that
$d_1\leq d_2 \leq \cdots \leq d_k$. Then using  Bertini Theorem, we can replace each equation $f_j$ by a linear combination
$$f_j'=f_j+\sum _{i=1,j-1}a_if_i$$
where $a_i$ is a generic homogeneous polynomial of degree $d_j-d_i$, such that if we define
$$Z_j=V(f_j',f_{j+1}',...,f_k')$$
then $Z_j$ is a complete intersection of dimension $(n+j-1)$ whose singular locus is contained in the finite set
$V_{sing}$, for all $j=1,2, \cdots, k$. To keep the notation simple, we assume in the sequel that the original polynomials $f_1, ..., f_k$ already satisfy this property.

\subsection{Step 2: non-proper Morse Theory and plan of proof} \label{s2b}

The main technical tool in proving our claim is the following result of Hamm,
see Proposition 3 in \cite{H2} for a more general version,  with our addition
concerning
the condition $(c0)$ in  Lemma 3 and Example 2 in \cite{DPtop}. See also  Proposition 11 in 
\cite{DP}.

\begin{prop} \label{prop1}

Let $A$ be a locally complete intersection in $\C^p$ with $\dim A=e$.
Let $g_1,...,g_p $ be polynomials in $\C[X_1,...,X_p]$ such that
the singular locus of $A$ is contained in $A_1=\{z \in A \ | \ g_1(z)=0\}$.
For $ 1 \leq j \leq p$, denote by  $\Sigma _j$ the set of critical
points
of the mapping $(g_1,...,g_j): A \setminus A_1 \to
\C^j$
and let $\overline \Sigma _j$ denote the closure of $\Sigma _j$ in $A$.
Assume  that the following conditions hold.

\begin{enumerate}

\item (c0) The set $\{z \in A \ | \ |g_1(z)| \leq a_1, ..., |g_p(z)| \leq a_p \}$
is compact for any positive real numbers $a_j$, $j=1,...,p$.

\item (c1) The  critical set $\Sigma _1$ is finite.

\item (cm) (for $m=2,...,p$) The map  $(g_1,...,g_{m-1}): \overline \Sigma _m  \to \C
^{m-1}$
is  proper.

\end{enumerate}

\noindent Then $A$ has the homotopy type of a space obtained from
$A_1$ by attaching $e$-dimensional cells.

\end{prop}

The main difference of the use of this result here versus its uses in \cite{DPtop} and \cite{DP} is that the
former transversality (or genericity) assumptions are replaced with the assumptions in Theorem  \ref{thm1}.

\medskip

Our plan of proving  Theorem  \ref{thm1} is  the following. Let
$$F_j=\{x \in \C^{N+1}~ | ~ f_j(x)=f_{j+1}(x)= \cdots = f_k(x)=0, ~ h(x)=1 ~\}$$
with the convention $F_{k+1}=F_h$, the Milnor fiber of $h$. In the next step we show that each of the pairs
$(F_{j+1},F_j)$ are in the situation of the pair $(A,A_1)$ in the above Proposition.
It follows that $F_{j+1}$ has the  homotopy type of a space obtained from $F_j$ by 
 attaching $(n+j)$-cells. Since $s \leq n-1 $, it follows that the $s$-connectivity of 
 $F_{j+1}$ implies the  $s$-connectivity of $F_j$, e.g. by using the homotopy exact sequence of the pair
$(F_{j+1},F_j)$. This proves the second claim in Theorem  \ref{thm1} modulo the following.

\subsection{Step 3: the study of the pairs $(F_{j+1},F_j)$ } \label{s2c}

As explained above, in this step we take $A=F_{j+1}=\{x \in \C^{N+1}~ | ~ f_{j+1}(x)= \cdots = f_k(x)=0, ~ h(x)=1 ~\}$,
$g_1=f_j$ and show that we can find generic linear forms $g_2,...,g_p$ with $p=N+1$
satisfying all the assumptions in Proposition  \ref{prop1}. More precisely, let $\SS_1$ be a Whitney regular stratification of the complete intersection $W_1=Z_j \cap H$, having as open stratum the regular part $S_0=W_{1,reg}$ of $W_1$. We choose 
$g_2$ such that the hyperplane $H_2:g_2=0$ is transverse to any stratum of the stratification $\SS_1$ and 
$\dim (V(g_1) \cap H_2)=N-2$. Then there is an induced Whitney regular stratification $\SS_2$ of $W_2=W_1\cap H_2$, whose strata are the non-empty intersections of the strata in $\SS_1$ with the hyperplane $H_2$. Choose $g_3$ such that the hyperplane $H_3:g_3=0$ is transverse to any stratum of the stratification $\SS_2$ and $\dim (V(g_1)   \cap H_2 \cap H_3)=N-3$ . Continue in the same way to define $g_4,...,g_p$.

The condition $(c0)$ in Proposition  \ref{prop1} is clearly satisfied, since 
$\dim ( V(g_1)  \cap H_2 \cap...\cap  H_p)=N-p=-1$ implies that the mapping
$(g_1,g_2,...,g_p):\C^p \to \C^p$ is finite, in particular it is proper.

Assume that the condition $(c1)$ fails. Then, exactly as in the proof in the previous section, we may find an irreducible curve $C \subset \Sigma _1$ and $g_1|C$ is constant, say equal to $c \in \C^*$. Note that $x \in  \Sigma _1$ implies that 
there is a relation of the following type
$$ d_xg_1=\lambda d_xh+\sum_{i=j+1,k}\lambda _id_xf_i.$$
Since $x \in A \setminus A_1$, we have $g_1(x) \ne 0$ and the Euler relation shows that $\lambda \ne 0$.
Therefore the above relation can be re-written as
$$d_xh=\sum_{i=j,k}\lambda _id_xf_i$$
with $\lambda _j \ne 0$.

Let $[C]=\{[x] \in \PP^n~|~x \in C\}$ and assume that $[x_n]$ is a sequence of points in  $[C]$ converging to a point
$[y] \in {\overline {[C]}}$. Since $x_n \in C$, it follows that $h(x_n)=1$, and hence the sequence $x_n$ is bounded away from the origin.
Two cases are possible. The first case, when the sequence of norms $|x_n|$ is bounded, leads to a convergent subsequence
$x_{n_m}$ converging to a point $y \in {\overline C} \subset A$. Using the Euler relation we get $d= \lambda _j(y)\cdot d_j  \cdot c$, i.e. $ \lambda _j$ is regular and non-zero at $[y]$.

The second case, when  the sequence of norms $|x_n|$ is unbounded, leads to a subsequence
$y_m=x_{n_m}$ such that $|y_m|\to +\infty$ and the sequence $z_m=\frac{y_m}{|y_m|}$ converges to some point $z$.
It follows that $g_1(z)= \lim f_j(y_m)\cdot |y_m|^{-d_j}= \lim c\cdot |y_m|^{-d_j}=0$ and similarly $h(z)=0$.
It follows that $[z] \notin Z_{j,sing}$. In particular $d_zf_j\ne 0$, and hence $\lambda_j$ does not have a pole at
$[z]=[y]$.

In conclusion, we get exactly as in the previous section, a regular section $\lambda_j$ without any zeroes of the line budle $\OO (d-d_j)$ over the irreducible projective curve $ {\overline {[C]}}$, which is a contradiction.

Finally, assume that $m\geq 2$ and that the condition $(c(m-1))$ holds but the condition $(cm)$ fails, i.e. the map
$$(g_1,...,g_{m-1}): \overline \Sigma _m  \to \C^{m-1}$$
is  not proper. This means that there is a sequence of points $p_n \in \overline \Sigma _m $ such that
$|p_n| \to +\infty$ and $g_a(p_n) \to c_a $ for $1\leq a \leq m-1.$ Since $ \Sigma _m $ is dense in 
$ \overline \Sigma _m $, we may assume that $p_n \in  \Sigma _m$. Moreover $ \Sigma _ {m-1}$ is clearly contained in
$\Sigma _m $, and $p_n \in \Sigma _m \setminus  \Sigma _ {m-1}$ (at least for $n$ large enough), since the condition
$(c(m-1))$ holds.
It follows that the differentials
$$d_{p_n}f_{j+1}, ...,d_{p_n}f_k, d_{p_n}h,
d_{p_n}g_1,...,d_{p_n}g_{m-1}$$
are linearly independent and $d_{p_n}g_m$ is a linear combination of them.
Let $q_n=\frac{p_n}{|p_n|}$ and assume that this sequence converges to some point $q$. Then exactly as above we get
$[q] \in  W_1  \cap H_2 \cap...\cap  H_{m-1}$, $[p_n]=[q_n]  \in  W_{1,reg}  \cap H_2 \cap...\cap  H_{m-1}$
and $\lim [q_n]=[q]$. It follows, using once again the Euler relation, that 
$[q] \in  W_1  \cap H_2 \cap...\cap  H_{m-1} \cap H_m$ and the hyperplane $H_m$ contains the limit of tangent spaces
$T=\lim T_{[q_n]}(W_{1,reg}  \cap H_2 \cap...\cap  H_{m-1})$. By the Whitney (a)-regularity condition,
this limit $T$ contains the tangent space $T'$ at $[q]$ to the stratum in the stratification $\SS_{m-1}$ of
$ W_1 \cap H_2 \cap...\cap  H_{m-1}$ that contains the point $[q]$. Therefore $T' \subset  H_m$, which contradicts the fact that $H_m$ is transverse to the stratification $\SS_{m-1}$.
This completes the verification of the assumptions in Proposition  \ref{prop1}
and hence the proof of Theorem  \ref{thm1}.

\end{document}